\documentclass[xypic,amscd,syntonly,amssymb]{amsart}

\usepackage{graphicx}
\usepackage[all]{xy}
\usepackage{amssymb}
\usepackage{amsmath}

\usepackage{epsfig}
\theoremstyle{plain}
\newtheorem{Thm}{Theorem}
\newtheorem{Prop}[Thm]{Proposition}
\newtheorem{Cor}[Thm]{Corollary}
\newtheorem{Lem}[Thm]{Lemma}

 \theoremstyle{definition}
\theoremstyle{remark}

\numberwithin{equation}{section}

\begin{document}
 \title{ Action Integrals and discrete series}

 \author{ ANDR\'{E}S   VI\~{N}A}
\address{Departamento de F\'{i}sica. Universidad de Oviedo.   Avda Calvo
 Sotelo.     33007 Oviedo. Spain. }
 \email{vina@uniovi.es}
\thanks{This work has been partially supported by Ministerio de Educaci\'on y
 Ciencia, grant   MAT2007-65097-C02-02}
  \keywords{Orbit method, geometric quantization, coadjoint orbits, discrete series}

 \maketitle
\begin{abstract}
Let $G$ be a complex semisimple   Lie group and ${G}_{\mathbb R}$
a real form that contains a compact Cartan subgroup $T_{\mathbb
R}$. Let $\pi$ be a discrete series representation of $G_{\mathbb
R}$. We present geometric interpretations in terms of concepts
associated with the manifold $M:=G_{\mathbb R}/T_{\mathbb R}$ of
the constant $\pi(g)$, for $g\in Z(G_{\mathbb R})$. For  some
relevant particular cases, we prove that this constant  is  the
action integral around a loop of Hamiltonian diffeomorphims of
$M$. As a consequence of these interpretations, we deduce lower
bounds for the cardinal of the fundamental group of some subgroups
of ${\rm Diff}(M)$. We also geometrically interpret the values of
the infinitesimal character of the
 differential representation of $\pi$.

\end{abstract}
   \smallskip


  MSC 2000: Primary: 53D50, \; Secondary: 22E45


\section{Introduction}

An irreducible unitary representation of a Lie group is a discrete
series  representation if it can be realized as a direct summand
of the left regular representation \cite{kn0,B-S}. If a group
possesses discrete series representations, then it contains a
compact Cartan subgroup. Kostant and Langlands conjectured the
realization of the discrete series by $L^2$-cohomology  of
holomorphic line bundles over the quotient of the group by the
compact Cartan subgroup. This conjecture has been proved for
``most" discrete series representations by Schmid in \cite{Sch1}
and fully in \cite{Sch2}.

 On the other hand, the Orbit Method \cite{Kir,Vo2} suggests the
 existence of a correspondence between the space of coadjoint
 orbits of a Lie group and its unitary dual (i.e. the set of
 equivalence classes of unitary irreducible representations), and the existence of
 some relations between geometric
 properties of the orbit and properties of the representation.
In the spirit of the  Orbit Method  and using the geometric
construction of Schmid, we describe here  interpretations of some
invariants of discrete series representations in terms of
geometric concepts relative to the  corresponding orbits.

The correspondence between coadjoint orbits and the unitary dual
is bijective for connected simply connected nilpotent groups; but,
  in general, there is not such a bijection.  However, it is
possible to associate an irreducible unitary representation   to
each {\em hyperbolic} orbit of a reductive group (see
\cite{Vo1,Vo2}). In \cite{V1}, we started with a hyperbolic orbit
and then we analyzed the geometric meaning of some invariants of
the corresponding representation.
 Here, we  start with a discrete series representation of a
 semisimple Lie group and using its realization by
 $L^2$-cohomology we will give the geometric interpretations above
 mentioned, which allow us to obtain results about  the homotopy of some groups of
 diffeomorphisms.

To formulate our results, we review briefly the geometric
construction of Schmid  and   introduce  notations that will be
used in the sequel.


\smallskip

{\bf The Schmid construction}

Let $G$ be a complex connected semisimple Lie group and
${\mathfrak g}:={\rm Lie}(G)$.
 For $A\in{\mathfrak g}$ and $g\in G$, we put $g\cdot A:={\rm
Ad}_g(A)$. For $\xi\in {\mathfrak g}^*$,  let $g\cdot\xi:={\rm
Ad}^*_g(\xi)$.  Let ${\mathfrak g}_{\mathbb R}\subset {\mathfrak
g}$ be a real form (non necessarily compact). By $G_{\mathbb R}$
we denote the connected subgroup of $G$ with Lie algebra
${\mathfrak g}_{\mathbb R}$.
 We will assume that $G_{\mathbb R}$
contains a {\em compact} Cartan subgroup $T_{\mathbb R}$ with Lie
algebra ${\mathfrak t}_{\mathbb R}$. We denote by $K_{\mathbb R}$
a maximal compact subgroup of
  $G_{\mathbb R}$ such that $T_{\mathbb
R}\subset K_{\mathbb R}$.

Let $ \Delta$ be a positive root system of ${\mathfrak
t}:={\mathfrak t}_{\mathbb R}\otimes_{\mathbb R} {\mathbb C}$ in
${\mathfrak g}$. From now on,
   $\rho$ stands for the   half the sum of the positive roots, $\rho=(1/2)\sum_{\nu\in \Delta}\nu$. We set
$${\mathfrak u}=\bigoplus_{\nu\in{\Delta}}{\mathfrak
 g}^{-\nu},\;\;\;\; \bar{{\mathfrak u}}=\bigoplus_{\nu\in{\Delta}}{\mathfrak
 g}^{\nu},$$
 ${\mathfrak g}^{\nu}$ being the root space of $\nu$.
  A root $\nu$ is said to be compact if ${\mathfrak
  g}^{\nu}\subset{\mathfrak k}$, and noncompact otherwise.
 We denote by ${\mathfrak b}$ the Borel subalgebra ${\mathfrak t}\oplus{\mathfrak
 u}$.
 The flag variety of ${\mathfrak g}$ will be denoted by  ${\mathcal
B}({\mathfrak g})$, or simply by ${\mathcal B}$. It can be
identified with
  the complex smooth variety   $G/B$, where
 $B=N_G({\mathfrak b})$, the normalizer of ${\mathfrak b}$ in $G$.

 Let $\phi$ be an element of the weight lattice of ${\mathfrak
t}_{\mathbb R}$. The corresponding character of $T_{\mathbb R}$ is
denoted by $\Phi$. By complex linearity, $\phi$ can be extended to
a map from ${\mathfrak t}$ to ${\mathbb C}$ and finally to an
element of ${\mathfrak g}^*$ by putting $\phi_{|{\mathfrak
u}\oplus\bar{\mathfrak u}}=0$.   The character $\Phi$ extends to a
group homomorphism on
 the Borel subgroup $B$ as well.

 We define the following holomorphic line bundle over
$\mathcal{B}\simeq G/B$.
 $${\mathcal V}=G\times_B {\mathbb C} =\{ (g,\,z)\,|\, g\in G,\, z\in
{\mathbb C} \}/\sim,$$
 where $(g,\,z)\sim (gb,\,\Phi(b^{-1})z)$, with $b\in
 B$. The natural $G$-action on ${\mathcal V}$ covers the action of
 $G$ on $\mathcal{B}$; that is, ${\mathcal V}$ is a $G$-equivariant
 holomorphic bundle.

 The set of orbits of the   $G_{\mathbb R}$-action   on ${\mathcal B}$ is finite. Since the Cartan
 subgroup $T_{\mathbb R}$ is compact  the  $G_{\mathbb R}$-orbit
 of ${\mathfrak b}$ is an open orbit, which can be identify with
 the quotient $M:=G_{\mathbb R}/T_{\mathbb R}$. The complex structure of ${\mathcal B}$ gives rise to
 a  $G_{\mathbb R}$-invariant complex structure on $M$, with
 $T^{(1,0)}_{eT_{\mathbb R}}=\bar {\mathfrak u}$. Thus, $M$ is a
 homogeneous complex manifold acted transitively by the group
 $G_{\mathbb R}$.

 The restriction of the holomorphic line ${\mathcal V}$ to $M$ will be also denoted by
 ${\mathcal V}$. Denoting by $L_A$ the left
invariant vector field on $G_{\mathbb R}$ defined by $A\in
{\mathfrak g}_{\mathbb R}$, the space
 of the
 holomorphic sections of ${\mathcal V}$ can be identified with the
 space of the functions $f\in C^{\infty}(G_{\mathbb R})$ such
 that,
 $f(gt)=\Phi(t^{-1})f(g)$ for all $t\in T_{\mathbb R}$ (i.e., they
 are $\Phi$-equivariant) and satisfy  $L_C(f)=0$ for all $C\in {\mathfrak
 u}$. Obviously,
  the
$G_{\mathbb R}$-action on $M$
  lifts to the sheaf
$\mathcal{O}(\mathcal V)$ of germs of holomorphic sections of
${\mathcal V}$.

On the space of compactly supported ${\mathcal V}$-valued
$(0,*)$-forms on $M$ we have the operator
 $\bar\partial$,
$$\bar\partial:{\mathcal A}^{0,*}({\mathcal V})\to {\mathcal
A}^{0,*+1}({\mathcal V}).$$
 Furthermore, the group $G_{\mathbb R}$ acts on the space ${\mathcal A}^{0,i}({\mathcal
 V})$ by translation and the action commutes with the operator
 $\bar\partial$.
 By means of $G_{\mathbb R}$-invariant Hermitian metrics on $M$ and
on ${\mathcal V}$, we define the operator $\bar\partial^*$, the
formal adjoint of $\bar\partial$. The space of
 square integrable, $C^{\infty}$, ${\mathcal V}$-valued
$(0,i)$-forms on $M$ which belong to
${\ker}(\bar\partial)\cap{\ker}(\bar\partial^*)$ is denoted by
$H^i_{(2)}(M,\,{\mathcal O}(\mathcal V))$.

Denoting by $(\,.\,,\,.\,)$ the bilinear form on ${\mathfrak t}^*$
induced by the Killing form, we put $q$ for the   integer
\begin{align}\label{q}
&\sharp\{ \;\nu\in \Delta\,|\,\nu\,\,{\rm
compact},\,(\phi+\rho,\,\nu)<0 \} \, + \\
 &\sharp\{\nu\in \Delta\,|\,\nu\,\,{\rm
noncompact},\,(\phi+\rho,\,\nu)>0 \} \notag.
\end{align}

According to the Langlands conjecture, if $\phi+\rho$ is regular,
the action of $G_{\mathbb R}$ on ${\mathcal
H}:=H^q_{(2)}(M,\,{\mathcal O}(\mathcal V))$ is an irreducible
unitary representation $\pi$, equivalent to a discrete series
representation (see  \cite{B-S} for the omitted details).

Note that every $\sigma\in{\mathcal H}$ is, in fact, a
$\bar\partial$-closed smooth Dolbeault form, but ${\mathcal H}$
can not be identified with $H^q(M,\,{\mathcal O}({\mathcal V}))$,
since $M$ is not compact, in general.

\smallskip


{\bf Statement of main results}

Throughout, we will assume that the element $\phi+\rho\in
i{\mathfrak t}^*_{\mathbb R}$ is regular, and we will denote by
 ${\mathcal H}_{K_{\mathbb R}}$ the space of $K_{\mathbb R}$-finite vectors
\cite{kn0,Sch3} of ${\mathcal H}$. The differential
representation, on  ${\mathcal H}_{K_{\mathbb R}}$, of the above
irreducible unitary representation $\pi$ will be denoted by
$\pi'$.

The purpose of this paper is to give geometric interpretations of
some invariants of $\pi$ and $\pi'$. Using these results we will
find lower bounds for the cardinal of the homotopy groups of some
subgroups of $\rm{Diff}(M)$, the group of diffeomorphisms of $M$.

On $W:={\mathbb C}\otimes (\bigwedge^q {\mathfrak
 u})^*$ we consider the representation $\Psi$ of $T_{\mathbb R}$,
 tensor product of
   $\Phi$ by the $q$-exterior product of ${\rm Ad}^*$. We denote by ${\mathcal
 P}$ the ${\rm GL}(W)$-principal bundle over $M$ determined by
 $\Psi$, and by ${\mathcal W}$ the associated vector bundle with fiber $W$ (defined by the standard representation
 of  ${\rm GL}(W)$). Thus, the vector space ${\mathcal H}$ is
 contained in the space
 $\Gamma({\mathcal W})$ of sections of ${\mathcal W}$. These
 bundles will be the geometric framework for our developments.

\smallskip

 {\bf 1. The representation $\pi'$ and the ``quantization operators".}
 Using the root structure of ${\mathfrak g}$ and   $\Psi'$, the
differential representation of $\Psi$, we will define a
$G_{\mathbb R}$-invariant connection $\Omega$ on
 ${\mathcal P}$. The corresponding covariant derivative on
 sections of ${\mathcal W}$ is denoted $\nabla$. By means of $\nabla$ and $\Psi$, for
 each $A\in\mathfrak{g}_{\mathbb R}$, we define a first order differential operator ${\mathcal Q}_A$ which acts
 on the sections of $\mathcal W$. The first order term of ${\mathcal Q}_A$
 is $-\nabla_{X_A}$, with $X_A$ the vector field on $M$ defined by
 $A$, and the zeroth order term is related with the value  of $\Psi'$
 at $A$. Abusing of language, the operators $\mathcal{Q}_A$    will be
  called ``quantization operators" (see paragraph below Theorem \ref{piprime(A)}).
 We will prove in  {\em Theorem \ref{piprime(A)}} that the operators
 ${\mathcal Q}_A$ restricted to ${\mathcal H}_{K_{\mathbb R}}$
  are a representation of
 ${\mathfrak g}_{\mathbb R}$ equivalent to $\pi'$.

The universal enveloping algebra of ${\mathfrak g}$ is denoted by
${\rm U}({\mathfrak g})$ and its center by ${\mathcal Z}(\mathfrak
g)$. Let $\chi$ denote the infinitesimal character \cite{Kn} of
the ${\rm U}({\mathfrak g})$-module ${\mathcal H}_{K_{\mathbb
R}}$. Let $\{C_1,\dots,C_r\}$ be a basis of ${\mathfrak
t}_{\mathbb R}$ and $E_{\nu}$ a basis of ${\mathfrak g}^{\nu}$.
 According to Proposition 5.34 of \cite{Kn},
each element   $J\in {\mathcal Z}(\mathfrak g)$   is a linear
combination, $p(E_{-\nu},C_i,E_{\nu})$,  of terms of the following
two types:

 (a) $C^{m_1}_{i_1}\dots C^{m_r}_{i_r}$, with $m_j\in{\mathbb Z}_{\geq 0}$.

 (b) $E^{p_1}_{-\nu_{i_1}}\dots E^{p_l}_{-\nu_{i_l}}
   C^{n_1}_{i_{1}}\dots C^{n_r}_{i_r}  E^{q_1}_{\nu_{i_{1}}}\dots
   E^{q_l}_{\nu_{i_{l}}}E_{\nu},$ with $\nu_i,\nu\in\Delta$ and $q_j,n_j,q_j\in{\mathbb Z}_{\geq 0}$.

 We will prove that the corresponding differential operator
$p({\mathcal Q}_{E_{-\nu}},\,{\mathcal Q}_{C_i},\, {\mathcal
 Q}_{E_{\nu}})$ on  the space ${\mathcal
 H}$ is the scalar one defined by the constant $\chi(J)$
  ({\em Theorem \ref{Propchi}}).


  \smallskip

{\bf 2. Invariants defined by Schur's Lemma.} Given an element
$g_1$ in the center of $G_{\mathbb R}$, by Schur's Lemma
\cite{kn0}
 \begin{equation}\label{Schur}
 \pi(g_1)=\kappa\,{\rm Id}_{\mathcal H},
 \end{equation}
$\kappa$ being a complex number of modulus $1$. In this paper, we
will give geometric interpretations of the invariant $\kappa$
associated with $\pi$:

 (i) in terms of  ``evolution equations" for elements in ${\mathcal
 H}$, equations
generated by families of ``quantization operators";

 (ii) as a gauge transformation of the fibre bundle ${\mathcal P}$
 that
 is the time-$1$ map of a flow in ${\mathcal P}$ that preserves
 the connection.

For the explanation of (i) and (ii) we need to introduce some
notations. Henceforth, $\{g_t\,|\, t\in[0,\,1] \}$ stands for a
smooth curve in $G_{\mathbb R}$ with the initial point at $e$. For
brevity, such a curve will be called a {\em path} in $G_{\mathbb
R}$. We denote by $\{A_t\}\subset{\mathfrak g}_{\mathbb R}$ the
corresponding velocity curve, that is,
\begin{equation}\label{volocity}
A_t=\Dot g_tg_t^{-1}.
 \end{equation}
  By means of the $G_{\mathbb R}$-action on $M$, the path  $g_t$ determines an isotopy of $M$, which will be denoted
  by $\{\varphi_t\}_{t\in[0,\,1]}$ in the sequel; that is,
  \begin{equation}\label{varphit}
  \varphi_t(gT_{\mathbb R})=g_tgT_{\mathbb R}.
  \end{equation}

Given a path $\{g_t\}$ in $G_{\mathbb R}$, we can consider the set
of sections $\sigma_t$ of $\mathcal{W}$ defined by the following
``evolution equations":
 \begin{equation}\label{dsigmatdt}
 \frac{{\rm d}\sigma_t}{{\rm d}t}={\mathcal Q}_{A_t}(\sigma_t),\; \;   \;
 \sigma_0=\sigma.
 \end{equation}
If $g_1$ belongs to $Z(G_{\mathbb R})$, the center of $G_{\mathbb
R}$, then we will prove
   in {\em Theorem \ref{Thmsigma}} that $\sigma_1=\kappa\sigma$, for any $\sigma\in{\mathcal H}_K$.

\smallskip

For each $A\in{\mathfrak g}_{\mathbb R}$, the vector field $X_A$
on $M$
 determines a vector field $U(X_A)$ on ${\mathcal P}$ such that the
 Lie derivative of $\Omega$ with respect this vector field
 vanishes. So, a path $\{g_t\}$ in $G_{\mathbb R}$  defines the time-dependent vector
  field $U(X_{A_{t}})$, which in turn determines a flow
${\bf H}_t$ on ${\mathcal P}$. If $g_1\in Z(G_{\mathbb R})$, we
will prove that   ${\bf H}_1$ is the gauge transformation ${\bf
H}_1(p)=p\kappa$ ({\em Theorem \ref{PropH1}}).

\smallskip

 When $G_{\mathbb R}$ is compact and    $\phi$ is a regular dominant weight,
the representation $\pi$ is the one provided by the Borel-Weil
Theorem. In this case  $M$ is the flag variety ${\mathcal
B}(\mathfrak g )$ of ${\mathfrak g}$, a compact manifold
diffeomorphic to the coadjoint orbit of $\phi\in{\mathfrak g}^*$.
On $M$ we will
  consider the symplectic Kirillov structure $\varpi$
\cite{Kir0}.
 When $g_1\in Z(G_{\mathbb R})$,
$\{\varphi_t\}$ is a loop in ${\rm Ham}(M,\,\varpi)$, the
Hamiltonian group
 of $(M,\, \varpi)$ \cite{Mc-S, lP01}. By studying equation
 (\ref{dsigmatdt}), we will find that $\kappa$ is
the exponential of the symplectic action \cite{We,V} around the
loop $\{\varphi_t\}$ (see {\em Proposition
\ref{SymplecticAction}}).


\smallskip

{\bf 3. Lower bounds for the cardinal of the homotopy group of
some subgroups of ${\rm Diff}(M)$.}
  Let $\mathfrak{X}(M)$ denote
the Lie algebra of all vector fields on $M$.
 We will consider  subalgebras  of
 $\mathfrak{X}(M)$, denoted by $\mathfrak{X}'$,
    such that each $Z\in\mathfrak{X}'$ determines a
 vector field $U$ on ${\mathcal P}$ satisfying ${\mathfrak
 L}_U\Omega=0$ (${\mathfrak L}$ = Lie derivative). The precise
 conditions imposed to the algebras $\mathfrak{X}'$  are detailed in the third paragraph of Section \ref{Sect.Invariance}.
 Let ${\mathcal G}$ be a connected Lie subgroup of ${\rm
 Diff}(M)$,
 which contains the isotopies associated with paths in $G_{\mathbb R}$
 and such that
 ${\rm Lie}({\mathcal G})$ is subalgebra of some algebra
 $\mathfrak{X}'$.
  Using the above interpretation of $\kappa$ as a gauge transformation
  which is the final point of a curve
 consisting  of automorphisms of ${\mathcal P}$,
  we will prove that:
 $$\sharp(\pi_1({\mathcal G}))\geq\sharp\{\Psi(g)\,|\,g\in Z(G_{\mathbb
 R})\},$$
(see {\em Theorem \ref{pi1G}}).


When $G_{\mathbb R}$ is compact and    $\phi$ is a regular
dominant weight, the above result adopts the following form: Let
 ${\mathcal G}$ be a connected subgroup of the symplectomorphism group of
$(M,\,\varpi)$ which contains $G_{\mathbb R}$, then
 $\sharp\{\Phi(g)\,|\, g\in Z(G_{\mathbb R})\}$ is a lower bound of the cardinal  of
 $\pi_1({\mathcal G})$ ({\em Theorem \ref{ThmpiHam}}).


 Theorem 6 of \cite{V2} gives a lower bound for $\pi_1({\rm
 Ham}({\mathcal O}))$, where ${\mathcal O}$ is a coadjoint orbit of $\rm {SU}(n+1)$ diffeomorphic to a (partial) flag
 manifold of $G=\rm {SL}(n+1,\;{\mathbb C})$ (i.e. diffeomorphic to the the quotient of $\rm{SL}(n+1,\,\mathbb{C})$
  by a parabolic subgroup). This lower bound
 is $\leq 2N$, where $N$ is the minimal Chern number of the
 manifold on spheres.
  For the flag
 variety   $\mathcal{B}(\mathfrak{sl}(n+1,\,{\mathbb C}))$ the minimal Chern
 number is $2$ (see \cite[page 117]{Mc-S2}). Hence, the lower bound given by
 the mentioned theorem applied to   this flag manifold is $\leq 4$, for
 any $n\geq 1$.
   In contrast with the result of \cite{V2},
 the lower bound given for $\sharp\pi_1({\rm
 Ham}(M,\,\varpi))$ in Theorem \ref{ThmpiHam} (when it is applied
 to $G_{\mathbb R}=\rm{SU}(n+1)$)
  depends on $n$
 and also on $\phi$. For example, this lower bound is
  equal to $n+1$, when $\Phi$ is injective on $Z(\rm{SU}(n+1))$.


\medskip

This article is organized as follows. In Section
\ref{SectPrincipal}, we define the principal bundle ${\mathcal
P}$,  the connection $\Omega$ on it and the vector bundle
${\mathcal W}$. We prove Theorem \ref{piprime(A)}, which gives the
representation $\pi'$ in terms of the ``quantization operators".
We also prove Proposition \ref{SymplecticAction} that relates
$\kappa$ with the action integral around a Hamiltonian loop in
$M$. The final part of this section   concerns the interpretation
of the infinitesimal character in terms of   ``quantization
operators".

 In Section \ref{Sect.Invariance}, we prove Theorem \ref{PropH1}
 about $\kappa$ as a gauge transformation of ${\mathcal P}$. We
 also prove Theorems \ref{pi1G} and \ref{ThmpiHam} which give lower
 bounds for the cardinal of the homotopy groups of subgroups of ${\rm Diff}(M)$.


\section{The ``quantization operators"} \label{SectPrincipal}

Let $E_{\nu}$ be a basis of ${\mathfrak g}^{\nu}$, then the
decomposition
 ${\mathfrak g}={\mathfrak t}\oplus{\mathfrak u}\oplus\bar{\mathfrak
 u}$ gives rise to the following direct sum decomposition
 \begin{equation}\label{SumDirecDecom}
  {\mathfrak g}_{\mathbb R}={\mathfrak t}_{\mathbb R}\oplus {\mathfrak
 l},
 \end{equation}
  where
  \begin{equation}\label{frakl}
   {\mathfrak l}:=\bigoplus_{\nu\in\Delta}\Big({\mathbb
 R}(E_{\nu}+\bar{E_{\nu}})\Big)\oplus \bigoplus_{\nu\in\Delta}\Big(i{\mathbb
 R}(E_{\nu}-\bar{E_{\nu}})\Big).
  \end{equation}
 Here the bar designates conjugation with respect to the real form ${\mathfrak g}_{\mathbb R}$.
  The component of  $A\in {\mathfrak g}_{\mathbb R}$  in ${\mathfrak t}_{\mathbb
 R}$ will be denoted by $A_0$.
  The decomposition (\ref{SumDirecDecom})  is preserved by ${\rm
ad}(C)$, when $C\in{\mathfrak t}_{\mathbb
 R}$. Thus,
 \begin{equation}\label{tcdota0}(t\cdot A)_0=t\cdot
 A_0,\;\;\;\;\hbox{for
 all}\;\,t\in{T}_{\mathbb R}.
 \end {equation}

 \smallskip

 Since $\phi$ belongs to the lattice weight of ${\mathfrak t}_{\mathbb R}$, it gives rise to a group homomorphism
 $\Phi:T_{\mathbb R}\to {\rm U}(1)$. On the other
 hand, we have the adjoint representation ${\rm Ad}:T_{\mathbb R}\to
 {\rm Aut}({\mathfrak u})$ and the corresponding exterior product
 of its dual on $(\bigwedge^q {\mathfrak u})^*$, where $q$ is the nonnegative integer defined by (\ref{q}).
 On the vector space  $W={\mathbb C}\otimes (\bigwedge^q {\mathfrak
u})^*$ we will consider the representation of $T_{\mathbb R}$
tensor product $\Phi\otimes\big(\bigwedge^q\,{\rm
 Ad}\big)^*$. As we said,  this representation will be denoted   by $\Psi$.

For $A\in {\mathfrak g}_{\mathbb
 R}$, we define $h_A:G_{\mathbb R}\to\mathfrak{gl}(W)$ by the formula
 \begin{equation} \label{h_A}
 h_A(g)=\Psi'((g^{-1}\cdot A)_0),
  \end{equation}
where $\Psi'$ is the derivative of $\Psi$.
 By
(\ref{tcdota0}) $h_A$ induces a map on $M$, which will  also be
denoted by $h_A$.

It is straightforward to prove the following lemma:
 \begin{Lem}\label{LemPsi(t)h}
 For any $g\in G_{\mathbb
 R}$ and every $t\in T_{\mathbb R}$,
 $$\Psi(t)h_A(g)=h_A(g)\Psi(t).$$
 \end{Lem}

On the other hand,
 $$[h_A(g),\,h_C(g)]_{\mathfrak{gl}}=[\Psi'((g^{-1}\cdot
A)_0),\, \Psi'((g^{-1}\cdot C)_0)   ]_{\mathfrak{gl}}=\Psi'[
(g^{-1}\cdot A)_0,\, (g^{-1}\cdot C)_0],$$
  where
$[\,,\,]_{\mathfrak{gl}}$ is the bracket in the Lie algebra
$\mathfrak{gl}(W)$.  As ${\mathfrak t}_{\mathbb
 R}$ is an abelian Lie algebra, we have
  \begin{equation}\label{commhAhC}
[h_A,\,h_C]_{\mathfrak{gl}}=0.
  \end{equation}

Given $A\in {\mathfrak g}_{\mathbb
 R}$, we   set $R_A$ for right invariant vector field on $G_{\mathbb R}$ defined by $A$.
   The vector field on $M$ defined by $A$ will be denoted by $X_A$, as we said.
From the definitions, it follows that
 \begin{align}\label{casesbis}
 R_A(h_C)=&-h_{[A,\,C]}\\\label{cases}
X_A(h_C)=&-h_{[A,\,C]} \; \;(\hbox{$h_C$ considered as function on
$M$)}.
\end{align}

\smallskip

 We denote by ${\mathcal P}$ the following ${\rm GL}(W)$-principal bundle over $M=G_{\mathbb
 R}/T_{\mathbb R}$:
 $${\mathcal P}:=G_{\mathbb R}\times_{\Psi}{\rm GL}(W)\overset{{\rm pr}}{\longrightarrow}
 M.
  $$
That is, ${\mathcal P}=\{[g,\,\alpha]\,|\,g\in G_{\mathbb R}
,\,\alpha\in{\rm GL}(W)\}$ with
$[g,\,\alpha]=[gt,\,\Psi(t^{-1})\alpha]$, for $t\in T_{\mathbb
R}$.  The right ${\rm GL}(W)$-action on ${\mathcal P}$ will be
denoted
 by ${\mathcal R}$ and by $V_y$ the vertical vector field associated with
 $y\in\mathfrak{gl}(W)$.

 We denote by ${\mathcal L}$ the left natural $G_{\mathbb
 R}$-action on ${\mathcal P}$; this action gives to ${\mathcal P}$ the
  structure of $G_{\mathbb R}$-equivariant bundle. Given $A\in{\mathfrak g}_{\mathbb
 R}$, the  vector field on ${\mathcal P}$ determined by
 $A$ through ${\mathcal L}$ will be denoted by $Y_A$. So $({\mathcal
 L}_g)_*(Y_A)=Y_{g\cdot A}.$

\smallskip

For $C\in{\mathfrak t}_{\mathbb R}$,  the trivial curve $\{[g {\rm
e}^{\epsilon C},\, \Psi({\rm e}^{-\epsilon C})\alpha]
 \}_{\epsilon}$ in ${\mathcal P}$ defines the vector
 $$Y_{g\cdot C}[g,\,\alpha]-V_y[g,\,\alpha],$$
  with
 $y={\rm Ad}_{\alpha^{-1}}(\Psi'(C))$;  here
 ${\rm Ad}$ is adjoint action of ${\rm GL}(W)$ on
 its Lie algebra.
 Hence, the tangent space to $\mathcal{P}$ at $[g,\,\alpha]$ is
 \begin{equation}\label{tangentF}
 T_{[g,\,\alpha]}\mathcal{P} =
 \frac{    \{ Y_A[g,\,\alpha] \;|\; A\in{\mathfrak g}_{\mathbb R}
\}\oplus\{
 V_y[g,\,\alpha] \;|\;y\in \mathfrak{gl}(W) \}              }
  {   \{ Y_{g\cdot C}[g,\,\alpha]- V_{{\rm
Ad}_{\alpha^{-1}}(\Psi'(C))}[g,\,\alpha]\;|\; C\in{\mathfrak
 t}_{\mathbb R}  \}     }\,.
 \end{equation}

\smallskip

  We define  the following
  $\mathfrak{gl}(W)$-valued $1$-form on ${\mathcal P}$:
\begin{equation}\label{definOmega}
\Omega(Y_A[g,\,\alpha]+V_y[g,\,\alpha])= {\rm
Ad}_{\alpha^{-1}}\big(\Psi'((g^{-1}\cdot A)_0)\big)+y.
\end{equation}
Obviously,  $\Omega$ is  well-defined on the quotient
(\ref{tangentF}). It is  straightforward to check its invariance
under ${\mathcal L}$ and that ${\mathcal R}^*_{\alpha}\Omega={\rm
Ad}_{\alpha^{-1}}\circ\Omega$, for all $\alpha\in {\rm GL}(W)$.
 That is, we have the following proposition:
 \begin{Prop}\label{Proconnection}
 The $1$-form $\Omega$
  is a $G_{\mathbb R}$-invariant connection on the
  ${\rm GL}(W)$-principal bundle
 ${\mathcal P}$.
 \end{Prop}

We set ${\bf h}_A[g,\,\alpha]:={\rm Ad}_{\alpha^{-1}}h_A(g)$. By
Lemma \ref{LemPsi(t)h} together with (\ref{tcdota0}), it follows
that ${\bf h}_A$ is a $\mathfrak{gl}(W)$-valued well defined map
on ${\mathcal P}$. In this notation
\begin{equation}\label{Omega1}
\Omega(Y_A+V_y)={\bf h}_A+y.
\end{equation}

The horizontal lift of $X_A$
  at the point $[g,\,\alpha]$
   is denoted by $X_A^{\natural}[g,\,\alpha]$. From (\ref{Omega1}),
  it follows that
 \begin{equation}\label{HorizLift}
X_A^{\natural}[g,\,\alpha
]=Y_A[g,\,\alpha]+V_{y}[g,\,\alpha],\,\;\hbox{with}\;\, y=-{\bf
h}_{A}[g,\,\alpha].
 \end{equation}

The following lemma is immediate:
\begin{Lem}\label{LemhA}
Given $A, C\in{\mathfrak g}_{\mathbb R}$, $y\in\mathfrak{gl}(W)$
and $p$ a point of  ${\mathcal P}$,
$$Y_A({\bf h}_C)=-{\bf h}_{[A,\,C]}\;\;\;\hbox{and}\;\;\; V_y(p)({\bf
h}_A)=-[y,\,{\bf h}_A(p)]_{\mathfrak{gl}}.$$
 \end{Lem}

The curvature of the connection $\Omega$ will be denoted by ${\bf
K}$.  As the curvature is a tensorial $2$-form \cite{Kob-Nom}, the
following proposition determines ${\bf K}$:
\begin{Prop}\label{ProbfK}
 For all $A,C\in{\mathfrak g}_{\mathbb R}$,
\begin{equation}\label{Curvature}
 {\bf K} (Y_A,\,Y_C)=
  -{\bf h}_{[A,\,C]}
 \end{equation}

\end{Prop}

 {\it Proof.}
  By the structure equation
 $${\bf K}(Y_A,\,Y_C)=d\,\Omega(Y_A,\,Y_C)+[\Omega(Y_A),\,\Omega(Y_C)]_{{\mathfrak{gl}}} .$$
By (\ref{commhAhC}),
$[\Omega(Y_A),\,\Omega(Y_C)]_{{\mathfrak{gl}}}=0$.
 From    Lemma \ref{LemhA}, it follows that
 $$Y_A(\Omega(Y_C))= {\bf
h}_{[C,\,A]}=-Y_C(\Omega(Y_A)).$$
 On the other hand, $\Omega([Y_A,\,Y_C])=-\Omega(Y_{[A,\,C]})=-{\bf h}_{[C,\,A]}$.
Thus, (\ref{Curvature}) follows.

 \qed

\smallskip

We denote by $D$ the covariant derivative determined by $\Omega$.
From (\ref{HorizLift}) together  Lemma \ref{LemhA},
(\ref{commhAhC}) and Proposition \ref{ProbfK},   it follows
\begin{equation}\label{DhCYA}
D{\bf h}_C(Y_A)=-{\bf K}(Y_C,\,Y_A), \;\;\hbox{for all}\;\;
A,C\in{\mathfrak g}_{\mathbb R}.
\end{equation}

\smallskip

The    vector bundle on $M$ associated with ${\mathcal P}$,
$$G_{\mathbb
 R}\times_{\Psi} ({\mathbb
C}\otimes(\Lambda^q{\mathfrak u})^*\big),$$
 will be denoted
${\mathcal W}$. The elements of ${\mathcal W}$ are classes
$\langle g,\,w\rangle$, of pairs $(g,\,w)\in G_{\mathbb R}\times
W$, with $\langle g,\, w\rangle=\langle
gt,\,\Psi(t^{-1})w\rangle$, for all $t\in{T}_{\mathbb R}$.
 ${\mathcal W}$ is a $G_{\mathbb R}$-equivariant vector bundle
 with the $G_{\mathbb R}$-action $g'\cdot\langle g,\,w\rangle=\langle g'g,\,w\rangle$.

The $C^{\infty}$ sections $\sigma$ of ${\mathcal W}$
 can be identified with the $C^{\infty}$ $\Psi$-equivariant
functions
$s:G_{\mathbb R}\to W$,   i.e. $C^{\infty}$ functions that satisfy
 \begin{equation}\label{s(gt)}
 s(gt)=\Psi(t^{-1})s(g), \;\;\hbox{for all}\; g\in G_{\mathbb
 R}\;\;
\hbox{and}\;  t\in T_{\mathbb R}.
 \end{equation}
The section $\sigma$ is related with the $\Psi$-equivariant
function $s$ by the formula
 \begin{equation}\label{sigma-equiv}
\sigma(gT_{\mathbb R})=\langle g,\,s(g)\rangle.
 \end{equation}

The $G_{\mathbb R}$-action on the section $\sigma$ is given by
$(g\cdot\sigma)(x)=g\cdot\sigma(g^{-1}x)$, for all $x\in M$. From
(\ref{sigma-equiv}),  it follows   the following proposition,
which gives the action of the representation $\pi$ on
$\Psi$-equivariant functions:

  \begin{Prop}\label{PropInva}
   The $\Psi$-equivariant function associated with
$g\cdot\sigma$ is $s\circ L_{g^{-1}}$, where $L_{g^{-1}}$ is the
left multiplication by $g^{-1}$ in $G_{\mathbb R}.$  In
particular, if $\sigma\in{\mathcal H}$, then $\pi(g)(s)=s\circ
L_{g^{-1}}.$
 \end{Prop}

\smallskip

On the other hand, a $\Psi$-equivariant function $s$ determines a
pseudotensorial \cite{Kob-Nom} map $s^{\flat}:{\mathcal P}\to W$
of type standard; that is,
$s^{\flat}(p\beta)=\beta^{-1}s^{\flat}(p)$, for all $\beta\in{\rm
GL}(W)$. The maps $s$ and $s^{\flat}$ are related by
\begin{equation}\label{ssflat}
s^{\flat}[g,\,\alpha]=\alpha^{-1}s(g).
\end{equation}
 From (\ref{ssflat}), it follows
 \begin{equation}\label{YAsflat}
 Y_A(s^{\flat})=(R_A(s))^{\flat},\;\;\;\;\; V_{{\bf
 h}_A}(s^{\flat})=-(h_A(s))^{\flat},
 \end{equation}
where $V_{{\bf h}_A}(s^{\flat})$ is the map  $p\in{\mathcal
P}\longmapsto V_{{\bf h}_A(p)}(s^{\flat})\in W$
and $h_A(s)$ is the mapping on $G_{\mathbb R}$ which at $g$
 takes the value $h_A(g)(s(g))$. The function $h_A(s)$ is in fact $\Psi$-equivariant as a
consequence of Lemma \ref{LemPsi(t)h}.

The covariant derivative on sections of ${\mathcal W}$, induced by
the connection $\Omega$ will be denoted by $\nabla$.

\begin{Prop}\label{nablaXA}
Given $A\in{\mathfrak g}_{\mathbb R}$ and a section $\sigma$ of
${\mathcal W}$,   the $\Psi$-equivariant map associated with
$\nabla_{X_A}\sigma$ is $R_A(s)+h_A(s)$.
\end{Prop}
{\it Proof.} The pseudotensorial map associated with
$\nabla_{X_A}\sigma$ is $X_A^{\natural}(s^{\flat})$. From
(\ref{HorizLift})  together with (\ref{YAsflat}), it follows the
proposition.

 \qed

The endomorphism of ${\mathcal W}$ over the identity determined by
$h_A$ is denoted by $F_A$; that is,
 $F_A(\langle g,\,w\rangle)=\langle g,\,
 h_A(g)(w)\rangle$.

\begin{Prop}\label{CurvNabla}
For $A,C\in{\mathfrak g}_{\mathbb R}$  one has
 \begin{equation}\label{CurvProp}
 \big( [\nabla_{X_A},\, \nabla_{X_C}]-\nabla_{[X_A,\,X_C]}  \big)(\sigma)=-F_{[A,\,C]}(\sigma),
  \end{equation}
$\sigma$ being a section of ${\mathcal W}.$
\end{Prop}

{\it Proof.} Let $s$ be the equivariant function associated with
$\sigma$ and $\sigma_1:=\nabla_{X_C}\sigma$. By Proposition
\ref{nablaXA}, $s_1=R_C(s)+h_C(s)$. Analogously
$\nabla_{X_A}\sigma_1$ has as associated equivariant function to
$R_A(s_1)+h_A(s_1)$. Hence, the equivariant function defined
 by $\nabla_{X_A}(\nabla_{X_C}\sigma)$ is
$$R_A(R_C(s)+h_C(s))+h_A(R_C(s)+h_C(s)).$$
By (\ref{casesbis})
$$R_A(h_C(s))=-h_{[A,\,C]}(s)+h_C(R_A(s)).$$
So, $\nabla_{X_A}(\nabla_{X_C}\sigma)$ is associated with
 \begin{equation}\label{nablaAC}
 R_A(R_C(s))
 -h_{[A,\,C]}(s)+h_C(R_A(s))+ h_A(R_C(s))+h_A(h_C(s)).
 \end{equation}

The equivariant function determined by
$\nabla_{[X_A,\,X_C]}\sigma=-\nabla_{X_{[A,\,C]}}\sigma$
  is
 \begin{equation}\label{[XAAX]}
  -(R_{[A,\,C]}(s)
 +h_{[A,\,C]}(s)).
 \end{equation}

From (\ref{nablaAC}), (\ref{[XAAX]}) and (\ref{commhAhC}), it
follows that   the equivariant function
 corresponding to the left hand side of
 (\ref{CurvProp}) is $-h_{[A,C]}(s)$, and the proposition is proved.

 \qed

\smallskip

 By Proposition \ref{PropInva}, the representation
differential $\pi'$ of the discrete series representation $\pi$,
considered as an action on the space ${\mathcal E}$ of the
$\Psi$-equivariant functions which correspond to the elements of
${\mathcal H}_{K_{\mathbb R}}$, is given by
\begin{equation}\label{piprime}
\pi'(A)(s)=-R_A(s).
 \end{equation}
 From Proposition \ref{nablaXA} and taking into account that the
 section associated with $h_A(s)$ is $F_A(\sigma)$,
 we obtain the following theorem:
  \begin{Thm}\label{piprime(A)}
  For any $A\in{\mathfrak g}_{\mathbb R}$, $\pi'(A)$ acting on elements of ${\mathcal H}_{K_{\mathbb R}}$  is the
 operator
 \begin{equation}\label{PA}
{\mathcal Q}_A:=-\nabla_{X_A}+F_A.
 \end{equation}
\end{Thm}
 Given   a closed symplectic quantizable \cite{WO} manifold $N$,   each Hamiltonian vector field
 on $N$  has associated a
 differential operator called quantization operator \cite{Sn}, which acts on the sections of a line bundle
  over N  called a  prequantum bundle. On
 the other hand,
 if $G_{\mathbb R}$  is compact and $\phi$  is a regular dominant weight,
 then the integer $q$ is zero; so, ${\mathcal W}$ is a line bundle on $M$. In this
 case $M$ equipped with the the Kirillov
 symplectic structure $\varpi$ is quantizable and  ${\mathcal
 W}$ is a prequantum bundle.  Furthermore,
  the above operator ${\mathcal
 Q}_A$ is the quantization operator associated with the vector
 field $X_A$. Thus, we will also call ``quantization operators" the differential operators
 defined in (\ref{PA}) for the general case.

\smallskip

Let $\{g_t\}_{t\in[0,1]}$ be a path in $G_{\mathbb R}$. It
determines  the velocity curve $\{A_t\}$, defined in
(\ref{volocity}). This path also gives rise to the isotopy
$\{\varphi_t\}$ on $M$ defined by (\ref{varphit}). Moreover, the
time-dependent vector field on $M$ determined  by $\varphi_t$ is
$X_{A_t}$; that is,
$$\frac{{\rm d}\varphi_t}{{\rm d}t}=X_{A_t}\circ\varphi_t.$$

 We consider the evolution equation (\ref{dsigmatdt}) for sections of ${\mathcal
 W}$.
 By Proposition \ref{nablaXA}, the $\Psi$-equivariant function
$s_t$ associated with the solution of
 (\ref{dsigmatdt}) satisfies
 \begin{equation}\label{dstdt}
 \frac{{\rm d}s_t}{{\rm d}t}=-R_{A_t}(s_t) ,\; \;   \;
  s_0= s.
 \end{equation}
 It is easy to prove the following proposition:
\begin{Prop}\label{Propst} The solution of (\ref{dstdt}) is
$s_t=s\circ L_{g_t^{-1}}.$
\end{Prop}

By Proposition \ref{PropInva},  when $s\in{\mathcal E}$
Proposition
 \ref{Propst} can be rephrased as saying that the solution of
(\ref{dstdt}) is
\begin{equation}\label{rephrased}
s_t=\pi(g_t)(s).
 \end{equation}

 \begin{Prop}\label{kappa=Psi} If $g_1\in Z(G_{\mathbb R})$, and $s$
is the equivariant function associated with an element
$\sigma\in{\mathcal H}$, then $s_1=\kappa s$.
 \end{Prop}
 {\it Proof.}
 The proposition   follows from (\ref{rephrased}) together with (\ref{Schur}).

 \qed

\begin{Cor}\label{Psi(g1)=kappa}
If $g_1\in Z(G_{\mathbb R})$, then $\Psi(g_1)=\kappa\,{\rm Id}_W$.
\end{Cor}
{\it Proof.} By Proposition
 \ref{Propst}, we have for all $g\in G_{\mathbb R}$
 $$s_1(g)=s(g_1^{-1}g)=s(gg_1^{-1})=\Psi(g_1) s(g),$$
  and the corollary follows from Proposition \ref{kappa=Psi}.

\qed

By the equivalence between (\ref{dsigmatdt}) and (\ref{dstdt}),
Proposition  \ref{kappa=Psi} gives rise to the following theorem:
\begin{Thm}\label{Thmsigma}
 If $\{g_t\}_{t\in[0,\,1]}$ is a path in $G_{\mathbb R}$ with
 $g_1\in Z(G_{\mathbb R})$ and $\sigma\in{\mathcal H}$, then the
 solution of the evolution equation (\ref{dsigmatdt}) satisfies
 $\sigma_1=\kappa\sigma$.
 \end{Thm}

\smallskip

 Let $x$ be a point of  $M$, we set $x_t:=\varphi_t(x)$.
 Evaluating
 (\ref{dsigmatdt}) at the point $x_t$, and we obtain
 \begin{equation}\label{fracrmd}
  \frac{{\rm d}\sigma_t}{{\rm d}t}(x_t)=-\iota({X_{A_t}}(x_t))\nabla \sigma_t+F_{A_t}(\sigma_t(x_t)).
  \end{equation}

Let $\{\mu_a\}_a$  be a local frame for ${\mathcal W}$ on an open
set $U$,  the solution $\sigma_t$ of (\ref{dsigmatdt}) can be
written as
 \begin{equation}\label{sigmat}
  \sigma_t=\sum_a m^a_t\mu_a,
  \end{equation}
   with $m_t^a$ a complex function
 defined on $U$. We denote by $\vartheta$ the connection form
 in this frame;
  that is,
 $\nabla\mu_a=\sum_c\vartheta^c_a\mu_c$. Thus,
  \begin{equation}\label{tildenabla}
 \nabla \sigma_t=\sum_a\big( {\rm d}m_t^a+\sum_cm^c_t\vartheta^a_c
 \big)\mu_a.
  \end{equation}
On the other hand, $F_{A_t}(\mu_a)=\sum_c(F_t)^c_a\mu_c$, with
$(F_t)^c_a$ a complex function defined on $U$.

If $\{x_t\}$ is contained in $U$,  let $m^a(t):=m^a_t(x_t)$. It
follows from (\ref{fracrmd}), (\ref{tildenabla}) and
(\ref{sigmat})
\begin{equation}\label{dmat}
\frac{{\rm d}m^a(t)}{{\rm d}t}=
\sum_c\big(-\vartheta^a_c(X_{A_t}(x_t))+(F_{t})^a_c(x_t)\big)m^c(t).
\end{equation}

If ${\rm dim}\,W=1$, then the above expression reduces to
 \begin{equation}\label{dmatsimple}
\frac{{\rm d}m(t)}{{\rm d}t}= \big(-\vartheta(X_{A_t}(x_t))+
h_{A_t}(x_t)\big)m(t).
\end{equation}
 So,
$$m(t)=m(0)\,{\rm exp}\Big(\int_0^t\big(-\vartheta(X_{A_t}(x_t))+h_{A_t}(x_t)\big){\rm d}t
\Big).$$

\smallskip

Next, we consider the particular case when  $G_{\mathbb R}$ is
compact and $\phi$ is a regular dominant weight. In this case $M$
is the flag manifold of ${\mathfrak g}$ and it supports the
symplectic
  structure $\varpi$, given by the $2$-form
 \begin{equation}\label{Kirillov}
\varpi(X_A(gT_{\mathbb R}),\, X_C(gT_{\mathbb
R}))=\phi((g^{-1}\cdot[A,C])_0).
 \end{equation}
Now, as    $q$ equals zero,  ${\mathcal W}$ is a complex line
bundle and
 $${\bf h}_A([g,\,\alpha])=h_A(g)=\phi((g^{-1}\cdot A)_0).$$
 The
curvature ${\bf K}$ projects a $2$-form $\omega$ on $M$. By
Proposition \ref{ProbfK} (or by (\ref{CurvProp}))
$\omega(X_A,\,X_C)=-h_{[A,\,C]}$; thus, $\omega=-\varpi$. From
(\ref{DhCYA}), it follows that
\begin{equation}\label{dhC=}
 {\rm d}h_C=\varpi(X_C,\,.\,).
  \end{equation}
That is, $h_C$ is a
 Hamiltonian function associated with the vector field $X_C$.
  If $g_1\in Z(G_{\mathbb R})$, then the evaluation
curve $x_t$ is closed and nullhomologous in $M$ \cite{L-M-P,
Mc-S}. So, by Stokes theorem
 \begin{equation}\label{ActionIntegral}
 \frac{m(1)}{m(0)}={\rm exp}\Big(\int_S \varpi   + \int_0^1
h_{A_t}(x_t)\,{\rm d}t \Big),
 \end{equation}
$S$ being a $2$-chain whose boundary is $x_t$.

  The right hand side
of (\ref{ActionIntegral}) is the exponential of the action
integral around the loop ${\varphi}_t$ \cite{We,V}. From Theorem
\ref{Thmsigma}, we deduce the following proposition:
 \begin{Prop}\label{SymplecticAction}
If $G_{\mathbb R}$ is compact,   $\phi$ is a regular dominant
weight and $g_1\in Z(G_{\mathbb R})$, then
 $\kappa$ is the exponential of the action integral around
the Hamiltonian loop ${\varphi_t}$.
 \end{Prop}

\smallskip

From now on until  the end of this section, we will concern with
the infinitesimal character of the $U({\mathfrak g})$-module
${\mathcal H}_{K_{\mathbb R}}$.

The ``representation" of the associative algebra ${\rm
U}({\mathfrak g})$ induced by $\pi'$  on the space ${\mathcal
H}_{K_{\mathbb R}}$ will  also be denoted by $\pi'$.

We denote the extension of $A\in{\mathfrak g}_{\mathbb R}\mapsto
\Psi'(A_0)\in\mathfrak {gl}(W)$ to a map on ${\mathfrak g}$  by
$\psi$; that is,
\begin{equation}\label{psi}
 \psi(A+iB)=\Psi'(A_0)+i\Psi'(B_0).
\end{equation}

\begin{Prop}\label{Proppsi}
  If $C\in{\mathfrak t}$ and $s$ is a $\Psi$-equivariant
function of ${\mathcal E}$ (i.e.   associated with an element
  of ${\mathcal H}_{K_{\mathbb R}}$), then
 $\pi'(C)\,s=\psi(C)\,s.$
 \end{Prop}

{\it Proof.}   By (\ref{psi}), we can assume that $C$ is an
element of ${\mathfrak t}_{\mathbb R}$. The proposition follows
from (\ref{piprime}) together with the fact that $s$ is
$\Psi$-equivariant.

 \qed

  Since $\pi$ is unitary and irreducible, the Harish-Chandra module
 ${\mathcal H}_{K_{\mathbb R}}$ has an infinitesimal character
 $\chi:{\mathcal Z}({\mathfrak g})\to{\mathbb C}$ (see \cite[Corollary 8.14]{kn0}).

 Let $\{C_1,\dots, C_r\}$ be a basis of ${\mathfrak t}$, and
 $$J'=\sum z_{m_1\dots
m_r}^{i_1\dots i_r} \,C^{m_1}_{i_1}\dots C^{m_r}_{i_r}\in{\mathcal
Z}(\mathfrak g),$$ where $z_{m_1\dots m_r}^{i_1\dots i_r}$ is a
complex number and $m_j$ a positive integer. By Proposition
\ref{Proppsi}, $\pi'(J')$ acting on ${\mathcal E}$ is the operator
 \begin{equation}\label{ConstOperator}
 \sum z_{m_1\dots
m_r}^{i_1\dots i_r}\,\big(\psi(C_{i_1})\big)^{m_1}\dots
\big(\psi(C_{i_r})\big)^{m_r},
  \end{equation}
  which is the multiplication by the {\em constant} $\chi(J')$, since $J'$ belongs to ${\mathcal Z}({\mathfrak g})$.
  This result expressed in terms of the corresponding
  ``quantization operators" gives rise to the following proposition:

\begin{Prop}\label{PropQuantCharacter}
 If $J'=\sum z_{m_1\dots
m_r}^{i_1\dots i_r}\,C^{m_1}_{i_1}\dots C^{m_r}_{i_r}\in{\mathcal
Z}(\mathfrak g)$, then
 the  differential operator
  $$\sum z_{m_1\dots
m_r}^{i_1\dots i_r}\,({\mathcal Q}_{C_{i_1}})^{m_1}\dots
({\mathcal Q}_{C_{i_r}})^{m_r},$$
 acting on ${\mathcal H}_{K_{\mathbb R}}$ is  the scalar operator defined by the constant
 $\chi(J')$.
\end{Prop}

As ${\mathfrak t}$ is abelian, the algebra ${\rm U}({\mathfrak
t})$
 can be identified with the algebra of polynomial functions on ${\mathfrak
 t}^*$. In this way, the constant (\ref{ConstOperator}) is
 $J'(\psi)$, the evaluation at $\psi$ of the corresponding polynomial.

Next, we will extend the result stated in Proposition
\ref{PropQuantCharacter} to general elements of ${\mathcal
Z}({\mathfrak g})$.

$\rho$, the half the sum of the positive roots, induces the
 linear map $C\in{\mathfrak t}\mapsto C-\rho(C)\in{\rm U}({\mathfrak
 t})$. By the universal property of ${\rm U}({\mathfrak t})$, it extends to an algebra
 automorphism $\tau$ of ${\rm U}({\mathfrak t})$.

 By Proposition 5.34 in \cite{Kn},
  if $J\in{\mathcal Z}({\mathfrak g})$, then
  \begin{equation}\label{frakP}
  J\in {\rm
U}({\mathfrak t})\oplus {\mathfrak P}, \;\;\hbox{where}\;\;\;
 {\mathfrak P}=\sum_{\nu\in \Delta}{\rm U}({\mathfrak g} )E_{\nu}.
 \end{equation}
 The projection of $J$ into ${\rm U}({\mathfrak t})$ will be denoted $\tilde
 J$.

The algebra of the elements in ${\rm U}({\mathfrak t})$
 invariant under the Weyl group $N_G(T)/T$ (with $T=\,\hbox{complexification
 of}\;T_{\mathbb R}$)
 will be denoted by ${\rm U}({\mathfrak
 t})^I$.
The
 Harish-Chandra isomorphism $\gamma:{\mathcal Z}({\mathfrak
 g})\to {\rm U}({\mathfrak t})^I$ is defined by $\gamma(J)=\tau(\tilde J)$ \cite{Kn}.

   As the Harish-Chandra module ${\mathcal H}_{K_{\mathbb R}}$ has an infinitesimal
   character, there exists $\lambda\in{\mathfrak
   t}^*$ such that
    \begin{equation}\label{InfChar}
   \chi(J)=\gamma(J)(\lambda),
    \end{equation}
     for all $J\in{\mathcal Z}({\mathfrak g})$. In  particular, for the element $J'$ considered above,
   $J'(\psi)$ must be equal to $\gamma(J')(\lambda)$; that is,
   $$J'(\psi)= (J'-J'(\rho))(\lambda).$$
    Since $J'(\rho)$, the evaluation of the polynomial $J'$ on $\rho$,  is a complex number (i.e.
   a constant polynomial), $(J'-J'(\rho))(\lambda)=J'(\lambda)-J'(\rho).$
   Hence, $\lambda=\psi+\rho$ is the element of ${\mathfrak t}^*$    which determines
   the infinitesimal character $\chi$ by means of (\ref{InfChar}).

By (\ref{frakP}),  any   element $J\in{\mathcal Z}({\mathfrak g})$
is of the form
 \begin{equation}\label{J}
  J=\sum ({\rm const})\,C^{m_1}_{i_1}\dots C^{m_r}_{i_r}+
 \sum ({\rm const})\,E^{p_1}_{-\nu_{i_1}}\dots E^{p_l}_{-\nu_{i_l}}
   C^{n_1}_{i_1}\dots C^{n_r}_{i_r}  E^{q_1}_{\nu_{i_1}}\dots
   E^{q_l}_{\nu_{i_l}}E_{\nu}.
   \end{equation}
   For abbreviation, we write $J=p(E_{-\nu},C_i,E_{\nu}).$

As the space ${\mathcal H}_{K_{\mathbb R}}$ is dense in ${\mathcal
H}$ \cite{Sch3}, we can state the following theorem that relates
$\chi(J)$ with the corresponding composition of ``quantization
operators" acting   on ${\mathcal H}$:

\begin{Thm}\label{Propchi}
 If $J\in{\mathcal Z}({\mathfrak g})$ is of the form (\ref{J}), then differential operator
 $$p({\mathcal Q}_{E_{-\nu}},\,{\mathcal Q}_{C_i},\, {\mathcal
 Q}_{E_{\nu}})$$
 acting on the space ${\mathcal H}$ is the scalar operator defined by the constant $\chi(J)=\gamma(J)(\psi+\rho)$.
\end{Thm}


\section{Homotopy groups of some subgroups of ${\rm Diff}(M)$}\label{Sect.Invariance}

In this section, we assume that the spaces of smooth functions
between two manifolds are endowed with the Whitney
${C}^1$-topology; in particular, the Lie algebras of vector fields
on a manifold.

We will consider   subalgebras $\mathfrak{X}'$ of
$\mathfrak{X}(M)$, the Lie algebra consisting of the vector fields
on $M$, such that each of its elements admits a lift to a ${\rm
GL}(W)$-invariant vector field on ${\mathcal P}$ which is an
infinitesimal symmetry of the connection $\Omega$. We will deal
with the following points related with  such an algebra
${\mathfrak X}'$:
 (i) Let $\{\psi_t\,|\,t\in[0,\,1]\}$ be a loop in the group ${\rm Diff}(M)$
 generated by a vector field of $\mathfrak{X}'$; we prove that
 there is a lift of $\psi_t$ to a flow ${\bf H}_t$ in ${\mathcal
 P}$, such that ${\bf H}_1$ is a gauge transformation of ${\mathcal
 P}$ which preserves $\Omega$. (ii) We check that ${\bf H}_1$ is
 the multiplication by a constant on each orbit $\{{\bf
 H}_t(p)\}_t$. (iii) We will construct loops in ${\rm Diff}(M)$
 such that the corresponding  transformations ${\bf H}_1$ are the
 multiplication by a constant on the total space ${\mathcal P}$.
 (iv) If $\varphi,\varphi'$ are loops of those considered in (iii),
 such that the corresponding constant are different, then
 $\varphi$ and $\varphi'$ won't be homotopic in any Lie subgroup
 of ${\rm Diff}(M)$ whose Lie algebra is contained in
 ${\mathfrak X}'$.

\smallskip

We state the precise conditions imposed to the algebras
${\mathfrak X}'$.
 Let  $\mathfrak{X}'$ be  a Lie subalgebra of $\mathfrak{X}(M)$
 satisfying the following conditions:

 1) There is a continuous  ${\mathbb R}$-linear map
 $$Z\in \mathfrak{X}'\mapsto U(Z)\in \mathfrak{X}({\mathcal P})$$
 such that $Z^{\natural}$ is the horizontal part of $U(Z)$.

 2) For each  $Z\in \mathfrak{X}'$, there is a $C^{\infty}$ map
 $a(Z):{\mathcal P}\to\mathfrak{gl}(W)$ such that:

  2a) $a(Z)(p\beta)=\beta^{-1}a(Z)(p)\beta$, for all $\beta\in{\rm
  GL}(W),\, p\in{\mathcal P}.$ That is, $a(Z)$ is a
  pseudotensorial function on ${\mathcal P}$ of type ${\rm Ad}$
  \cite{Kob-Nom}.

  2b) $U(Z)=Z^{\natural}+V_{a(Z)}.$

  2c) $Da(Z)=-{\bf K}(Z^{\natural},\,.\,)$.

  2d) The map $Z\mapsto a(Z)$ is  ${\mathbb R}$-linear and continuous.

\smallskip

{\it Example.} For $X_C$, with $C\in{\mathfrak g}_{\mathbb R}$ we
define $U(X_C):=Y_C$ and $a(X_C):={\bf h}_C$. By (\ref{HorizLift})
and (\ref{DhCYA}), the Lie algebra ${\mathfrak X}':=\{X_C\,|\,
C\in {\mathfrak g}_{\mathbb R}\}$, with the above choices for
$U(X_C)$ and $a(X_C)$, satisfies the conditions 1)-2d).

\begin{Thm}\label{ThmLiederivative}
For any $Z\in{\mathfrak X}'$, the Lie derivative ${\mathfrak
L}_{U(Z)}\Omega$
 vanishes.
\end{Thm}

{\it Proof.} We will prove that the $1$-form ${\rm
d}(\iota_U\Omega)+\iota_U({\rm d}\Omega)$ is zero, where
$U:=U(Z)$. Let $p$ be an arbitrary point of ${\mathcal P}$, we
will check that
\begin{equation}\label{LieDerivativeE}
\big({\rm d}(\iota_U\Omega)+\iota_U({\rm d}\Omega)\big)(E)=0,
 \end{equation}
 when $E$ is a vertical vector at $p$ and when $E$ is horizontal.
 In any case, by 2b)
 \begin{equation}\label{anycase}
 {\rm d}(\iota_U\Omega)={\rm d}a,
 \end{equation}
 where $a:=a(Z)$.

 If $E$ is a vertical vector $E=V_y(p)$, with
 $y\in\mathfrak{gl}(W)$, by (\ref{anycase}) and 2a)
  $${\rm d}(\iota_U\Omega)(E)={\rm d} a(V_y)=-[y,\,a(p)].$$
 We extend $E$ to a vertical vector field. By the structure
 equation and 2b), we have
 $$(\iota_U({\rm d}\Omega))_p(E)={\bf K}_p(U,\,E)-[\Omega_p(U),\,\Omega_p(E)]=
 -[a(p),\,y].$$
 Thus, (\ref{LieDerivativeE}) holds when $E$ is a vertical vector.

 If $E$ is horizontal, then by (\ref{anycase}) and 2c)
 $${\rm d}(\iota_U\Omega)(E)=-{\bf K}(Z^{\natural},E).$$
 Now, we extend $E$ to an horizontal field; by the structure
 equation
 $$\big(\iota_U({\rm d}\Omega)  \big)(E)={\bf K}(U,\,E)={\bf
 K}(Z^{\natural},\,E).$$
 So, (\ref{LieDerivativeE}) also holds  when $E$ is a horizontal
 vector.

  \qed

Since $a(Z)$ satisfies 2a) and  $U(Z)=Z^{\natural}+V_{a(Z)}$, we
have the following proposition:
\begin{Prop}\label{Uinvariante}
For any $\beta\in{\rm GL}(W)$ and any $Z\in\mathfrak{X}'$,
$$({\mathcal R}_{\beta})_*(U(Z))=U(Z).$$
\end{Prop}

\smallskip

Let $\{Z_t\}_{t\in[0,\,1]}$ be a time-dependent vector field on
$M$, with $Z_t\in\mathfrak{X}'$. Let $\{\psi_t\}$ be the isotopy
defined by
\begin{equation}\label{isotopypsi}
\frac{{\rm d}\psi_t}{{\rm d}t}=Z_t\circ\psi_t,\;\;\;\;\psi_0={\rm
Id}_M.
\end{equation}
 As $M$ is a homogeneous space that admits a Hermitian metric, the flow
 $\psi$ is defined  on $[0,\,1]\times M$ \cite[page 85]{Lang}.

We have the time-dependent vector field $U_t:=U(Z_t)$ on
${\mathcal P}$ and the corresponding flow ${\bf H}$ which is also
defined on $[0,\,1]\times {\mathcal P}$:
\begin{equation}\label{bdH}
 \frac{{\rm d}{\bf H}_t(p)}{{\rm d}t}=U_t({\bf H}_t(p)),\;\;\; {\bf
 H}_0={\rm Id}_{\mathcal P}.
  \end{equation}
 From  Theorem \ref{ThmLiederivative}, it follows that the
diffeomorphism ${\bf H}_t$
 preserves the connection; that is, ${\bf H}_t^*\Omega=\Omega.$

 By 2b), given $p\in{\mathcal P}$, the curve ${\bf H}_t(p)$ satisfies ${\rm pr}({\bf
 H}_t(p))=\psi_t(x)$, if
 ${\rm pr}(p)=x$; that is, ${\bf H}_t(p)$ can be considered as   a lift of the curve $\psi_t(x)$ to ${\mathcal
 P}$ at the point $p$.
  In particular, if $\{\psi_t(x)\,|\,t\in[0,\,1]\}$ is a closed
 curve, then ${\bf H}_1(p)\in({\rm pr})^{-1}(x)$.

   Let us assume that  $\psi_0=\psi_1={\rm Id}_M$, that is
 $\{\psi_t\}_{t\in[0,\,1]}$ is a loop of diffeomorphisms of $M$;
 by Proposition \ref{Uinvariante},
${\bf H}_1$ is a gauge transformation of ${\mathcal P}$.
 Thus, there exists a map $f:{\mathcal P}\to{\rm GL}(W)$ satisfying
 $f(p\beta)=\beta^{-1}f(p)\beta$ and  such that
 ${\bf H}_1(p)=pf(p)$, for all $p\in{\mathcal P}$.
 As
 \begin{equation}\label{bfH1Omega}
  \Omega={\bf H}_1^*(\Omega)=f^{-1}{\rm d}f+{\rm Ad}(f^{-1})\circ
\Omega,
 \end{equation}
  we deduce
 \begin{equation}\label{dfHorizontal}
 {\rm d}f(E)=0,\;\;\;\hbox{if $E$ is horizontal.}
 \end{equation}
By Proposition \ref{Uinvariante}, $\Omega(U_t)={\rm
Ad}_{\beta^{-1}}(\Omega(U_t))$, for all $\beta\in{\rm GL}(W)$.
From (\ref{bfH1Omega}), it follows
$${\rm Ad}_{\beta^{-1}}(\Omega(U_t))= f^{-1}{\rm d}f(U_t)+{\rm Ad}(f^{-1})(\Omega(U_t)).$$
 Given $p\in{\mathcal P}$, taking $\beta=f(p)$, we conclude that
 \begin{equation}\label{dfUt}
 {\rm d}f(U_t)=0.
  \end{equation}
 Thus,   $f$ is {\em constant} on the orbit
  $\{{\bf H}_t(p)\,|\,t\in[0,\,1]\}$.

 \smallskip

Let $\{g_t\}$ be a path in  in $G_{\mathbb R}$ and $\{A_t\}$ the
corresponding velocity curve. Let us assume that
$X_{A_t}\in\mathfrak{X}'$.  This family gives rise to the
time-dependent vector field on ${\mathcal P}$, $U_t:=U(X_{A_t})$,
 which in turn defines a flow ${\bf
H}_t$   on ${\mathcal P}$. By (\ref{HorizLift}), $U(X_{A_t})
=Y_{A_t}$. So,
\begin{equation}\label{FlowH}
 \frac{{\rm d}{\bf H}_t}{{\rm d}t}=Y_{A_t}\circ {\bf H}_t,\;\;\;
 {\bf H}_0={\rm Id}_{\mathcal P}.
 \end{equation}
 \begin{Lem}\label{LemHt}
 The bundle diffeomorphism ${\bf H}_t$ defined in (\ref{FlowH}) is
the left multiplication by $g_t$ in ${\mathcal P}$; i.e. ${\bf
H}_t={\mathcal L}_{g_t}.$
\end{Lem}

{\it Proof.} Given $p\in{\mathcal P}$
 $$ \frac{{\rm d}}{{\rm
d}u}\bigg|_{u=t}g_up= \frac{{\rm d}}{{\rm
d}u}\bigg|_{u=t}g_ug_t^{-1}g_tp=Y_{A_t}(g_tp).$$

 \qed

 If $g_1$ is an element of the center of $G_{\mathbb R}$, then $\varphi_0=\varphi_1$ and  ${\bf
 H}_1$
  is a gauge transformation of ${\mathcal P}$.
\begin{Thm}\label{PropH1}
If $g_1\in Z(G_{\mathbb R})$, then for all $p\in{\mathcal P}$,
${\bf H}_1(p)=p\kappa$, where $\kappa$ is given by (\ref{Schur}).
 \end{Thm}
 {\it Proof.} By Lemma \ref{LemHt} together with Corollary \ref{Psi(g1)=kappa},
$${\bf
H}_1([g,\,\alpha])=[g_1g,\,\alpha]=[gg_1,\,\alpha]=[g,\,\Psi(g_1)\alpha]=[g,\,\alpha]\kappa.$$

 \qed

{\it Remark.}    Note that  paths with the same final point in
$Z(G_{\mathbb R})$ determine the same gauge transformation ${\bf
H}_1$.

\smallskip

Let ${\mathcal G}$ be a connected Lie subgroup of
 ${\rm Diff}(M)$  such that:
\begin{equation}
\begin{aligned}\label{aligned}
 &\hbox{(i)} \;{\rm Lie}( {\mathcal G})\subset\mathfrak{X}'. \\
 & \hbox{(ii)} \;\hbox{If} \;\{g_t\}_{t\in[0,\,1]}\; \hbox{is a path in}\; G_{\mathbb R},\;
 \hbox{then the isotopy}\;
\{\varphi_t\} \;
 \hbox{is contained in}\;
   {\mathcal G}.
\end{aligned}
\end{equation}
 We will  prove that
$$\sharp\pi_1({\mathcal G})\geq\sharp\{\Psi(g)\,|\,g\in Z(G_{\mathbb R})
\},$$
 but we need some previous results.

Let $\{g_t\}$ be a path in  in $G_{\mathbb R}$ such that
  $g_1\in Z(G_{\mathbb R})$.  Let
$\{\zeta^s\}$ be a {\em deformation} of the loop
$\varphi=\{\varphi_t\}$ in ${\mathcal G}$. That is, for each $s$,
$\zeta^s$ is a loop in ${\mathcal G}$ at ${\rm Id}_M$, with
$\zeta^0=\varphi$. We also assume that $s\mapsto \zeta^s$ is a
continuous map (considering   ${\mathcal G}$ equipped with the
$C^1$-topology, as we said).

For each $s$, we have the corresponding time-dependent vector
field $Z^s_t\in\mathfrak{X}'$, given by
$$\frac{{\rm d}\zeta^s_t}{{\rm d}t}=Z^s_t\circ\zeta^s_t.$$
The respective time-dependent vector field on ${\mathcal P}$,
$U^s_t:=(Z^s_t)^{\natural}+V_{a(Z^s_t)}$ determines the
corresponding flow ${\bf H}^s_t$. As above ${\bf H}^s_1$ is a
gauge transformation, which can be written ${\bf
H}^s_1(p)=pf^s(p)$.
By (\ref{dfHorizontal}),
 we have
\begin{equation}\label{df^s}
{\rm d}f^s(\hbox{Horizontal})=0.
 \end{equation}
 As before, $f^s$
   is constant on the orbit
 $\{{\bf
 H}^s_t(p)\,|\,t\in[0,\,1]\}$. Furthermore, this orbit is
 a lift to  ${\mathcal P}$ of curve
  $\{\zeta^s_t(x)\,|\,t\in[0,\,1]\}$, if $x={\rm pr}(p)$.

 In the statement of the following theorem we  consider the $T_{\mathbb R}$-principal bundle $T_{\mathbb
 R}\to G_{\mathbb R}\to M=G_{\mathbb R}/T_{\mathbb R}.$ We assume
 that this bundle is endowed  with the invariant connection
 \cite[page 103]{Kob-Nom}
 determined by the splitting (\ref{SumDirecDecom}).

 \begin{Thm}\label{ThmPerturb}
Let $\{g_t\}$ be a path in $G_{\mathbb R}$ which
 is a horizontal curve with
respect the invariant connection in $G_{\mathbb R}\to M$ and such
that its final point  belongs to $Z(G_{\mathbb R})$.  If
$\{\zeta^s\}_s$ is a deformation of
 $\varphi$ in ${\mathcal
G}$, then the gauge transformation ${\bf H}^s_1$ of ${\mathcal P}$
defined by $\zeta^s$ satisfies
 $${\bf H}^s_1(p)=p\kappa\;\; \hbox{for all}\;\; p\in({\rm pr})^{-1}(eT_{\mathbb R}),$$
 with $\kappa\,{\rm Id}_{W}=\Psi(g_1)$.

 \end{Thm}

{\it Proof.} Let $x_0$ denote  the  point $eT_{\mathbb R}\in M$.
We have the following closed curves in $M$:
 $$\{\varphi_t(x_0)\,|\,t\in[0,\,1]\},\;\;
 \{\zeta^s_t(x_0)\,|\,t\in[0,\,1]\}.$$
 By hypothesis, $\{g_t\}$ is the   horizontal lift
  of $\{\varphi_t(x_0)\}$  at the point $e\in G_{\mathbb R}$.

 Fixed $s$, for $a\leq s$, let $\gamma^a_t$ be the curve in $G_{\mathbb R}$  horizontal lift of $\{\zeta^a_t(x_0)\,|\,
  t\in[0,\,1]\}$ at the point $e$. We can construct a path $\{\Hat
  g_t\,|\, t\in[0,\,1]\}$ in
  $$\{\gamma^a_t\,|\,a\in[0,\,s],\;t\in[0,\,1]\}\subset G_{\mathbb R},$$
  with $\Hat g_1=g_1$ and
  such that the loop $\{\Hat
  g_t x_0\,|\, t\in[0,\,1]\}$ is as close to $\{\zeta^s_t(x_0)\,|\,
  t\in[0,\,1]\}$ as we wish; hence, $\{\Hat
  g_t x_0\}_t$ is contained in an arbitrarily small tubular neighborhood of
  $\{\zeta^s_t(x_0)\}_t$ in M (see Figure 1).

\begin{figure}[htbp]
\begin{center}
\epsfig{file=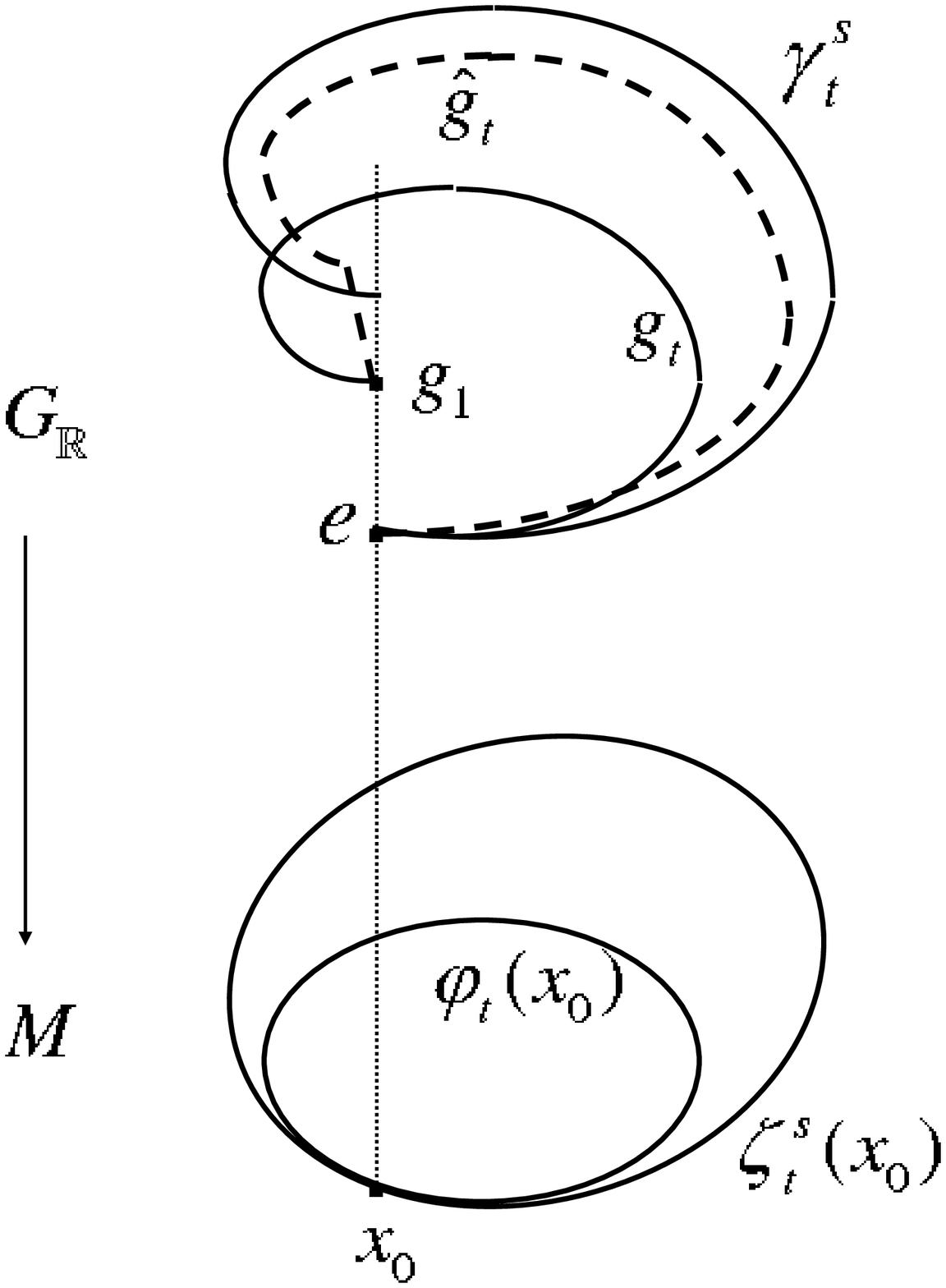, height=6.5cm}
 \caption[Figure 1]
 {\small The path $\Hat g_t$}
\end{center}
\end{figure}

  By (\ref{aligned})(ii), $\{\Hat g_t\}$ defines a loop in ${\mathcal G}$. With $\{\Hat g_t\}$
   we construct the corresponding flow $\Hat{\bf
   H}_t$ in ${\mathcal P}$ (see (\ref{FlowH})). Fix an arbitrary point $p$ in
    ${\mathcal P}$ such that
   ${\rm pr}(p)=x_0$.  The curves in ${\mathcal P}$ $\{\Hat{\bf
   H}_t(p)\}_t$ and
 $\{{\bf H}^s_t(p)\}_t$, which   are lifts at $p$ of
 $\{\Hat g_t x_0\,|\, t\in[0,\,1]\}$  and $\{\zeta^s_t(x_0)\,|\,
  t\in[0,\,1]\}$, will be  as   close to each other  as we will wish.

  Furthermore, $\Hat{\bf H}_1$ is defined by  a ${\rm
  GL}(W)$-valued map $\Hat f$ on ${\mathcal P}$. According to the Remark
  below Theorem \ref{PropH1}, $\Hat f=\kappa\,{\rm Id}_{{\rm
  GL}(W)}$, with $\kappa\,{\rm Id}_{W}=\Psi(g_1)$ (Corollary \ref{Psi(g1)=kappa}).

 Thus, for each tubular neighborhood of $\{{\bf H}^s_t(p)\,|\,
  t\in[0,\,1]\}$  in ${\mathcal P}$, there is a loop   in ${\mathcal G}$, defined by a path $\{\Hat g_t\}$,
  which in turn determines
  a flow $\{\Hat{\bf H}_t\}$ in ${\mathcal P}$ such that
  the corresponding evaluation curve $\{\Hat{\bf H}_t(p)\}$  is in that
  neighborhood and    the  gauge transformation $\Hat{\bf
  H}_1$ is the multiplication by $\kappa$. On the other hand, ${\bf H}_1^s$ on $\{{\bf
  H}^s_t(p)\,|\,t\in[0,\,1]\}$ is defined by the constant function
  $f^s$.
 By continuity,
  $f^s(p')=p'\kappa$, for any point $p'$ in $\{{\bf
  H}^s_t(p)\}_t$.


  \qed

 \begin{Cor}\label{CorPert} Let $\{g_t\}$ and $\{g'_t\}$ be
 paths in $G_{\mathbb R}$ satisfying the hypotheses of Theorem
 \ref{ThmPerturb}. If $\Psi(g_1)\ne\Psi(g'_1)$, then the
 corresponding loops $\varphi=\{\varphi_t\}$ and $\varphi'=\{\varphi'_t\}$ are not
 homotopic in ${\mathcal G}$.
 \end{Cor}

The Lie algebra of the holonomy  group at $e$ of the invariant
connection on $G_{\mathbb R}$ is generated by the  vectores of the
form $[A,\,C]_0$, with $A,C\in{\mathfrak l} $ (see \cite[Theorem
11.1]{Kob-Nom}). On the other hand, the  vectors $[E_{\nu},\,\bar
E_{\nu}]$, with $\nu\in\Delta$ span $i{\mathfrak t}_{\mathbb R}$.
So, from (\ref{frakl}), it follows that the
 Lie algebra of the mentioned holonomy  group is
${\mathfrak t}_{\mathbb R}$. Since $Z(G_{\mathbb R})\subset
T_{\mathbb R}$, each element of $Z(G_{\mathbb R})$ can be joined
to $e$ by a horizontal curve in $G_{\mathbb R}$ (horizontal with
respect to the invariant connection). Thus, from Corollary
\ref{CorPert}, it follows  the following theorem:
 \begin{Thm}\label{pi1G}
 If ${\mathcal G}$ is a connected Lie subgroup of ${\rm Diff}(M)$, for which
 the conditions (\ref{aligned})  hold,  then
 $$\sharp(\pi_1({\mathcal G}))\geq\sharp\{\Psi(g)\,|\,g\in Z(G_{\mathbb
 R})\}.$$
 \end{Thm}

\smallskip

As we said in Section \ref{SectPrincipal},   when $G_{\mathbb R}$
is compact and $\phi$ is a regular dominant weight,
$(M,\varpi)$ is the coadjoint orbit associated with $\phi$; thus,
  it is a     compact  symplectic manifold. In this case, as we
will see, a possible subalgebra ${\mathfrak X}'$ is the one
consisting of all locally Hamiltonian vector fields; that is, the
set of all vector fields $Z$ on $M$ such that ${\mathfrak
L}_Z\varpi=0$.

In fact, $M$ is simply connected \cite[page 33]{G-L-S}; so, such a
vector field $Z$ is Hamiltonian.  We denote by $a(Z)$ the
corresponding normalized Hamiltonian function; that is, $a(Z)$ is
the function on $M$ such that
$${\rm d}a(Z)=\iota_Z\varpi,\;\;\;\hbox{and}\;\;\int_Ma(Z)\varpi^{n}=0,$$
$2n$ being the dimension of $M$. We put $U(Z)$ for the vector
field on ${\mathcal P}$ defined by $U(Z)=Z^{\natural}+V_{a(Z)}.$
 As ${\bf K}$ projects on $M$ the form $-\varpi$ (see
 paragraph before (\ref{dhC=})), it is easy to check that
  the algebra $\mathfrak{X}'$ of Hamiltonian vector fields on
$M$, the Hamiltonian functions $a(Z)$ and the vector fields $U(Z)$
satisfy the conditions 1), 2a) - 2d) introduced at the beginning
of this Section.

Since any   Lie subgroup ${\mathcal G}$ of ${\rm
Symp}(M,\,\varpi)$, the group of all symplectomorphisms of
$(M,\,\varpi)$,
  satisfies the condition (\ref{aligned})(i), it follows from
  Theorem \ref{pi1G} the following result:
  \begin{Thm}\label{ThmpiHam}
Let us assume that $G_{\mathbb R}$ is compact and that $\phi$ is a
regular dominant weight. If   ${\mathcal G}$ is  any connected Lie
subgroup of ${\rm Symp}(M,\,\varpi)$ which contains $G_{\mathbb
R}$, then
 $$\sharp(\pi_1({\mathcal G}))\geq\sharp\{\Phi(g)\,|\,g\in Z(G_{\mathbb
 R})\}.$$
  \end{Thm}

\smallskip

{\it Remark.}  For $G_{\mathbb R}={\rm SU}(2)$, the corresponding
flag manifold is ${\mathbb C}P^1$. Given $[z_0:z_1]\in {\mathbb
C}P^1$ with $z_0\ne 0$,   we put $x+iy=z_1/z_0$. It is easy to
check that the vector fields $X_C$ and $X_D$ defined by the
matrices of $\mathfrak{su}(2)$
$$C=\,\rm{skew \;diagonal}\,(-c,\,c),\;\; D=\,\rm{skew\;
diagonal}\,(di,\,di),$$
 take at the point $(x=0,\,y=0)$ the values
 $X_C=-c\,\partial_x$, $X_D=d\,\partial_y.$ We denote by $\omega$
 the Fubini-Study form on ${\mathbb C}P^1$, then
 $\omega_{[1:\,0]}(X_C,\,X_D)=-cd/\pi$.

 On he other hand,
  let $\phi$ be the weight defined by
  $\phi({\rm diagonal}(ai,\,-ai))=a$. Then
  $\varpi_{[1:\,0]}(X_C,\,X_D)=\phi([C,D])=2cd$.  By the invariance
  of $\omega$ and $\varpi$ under the action of ${\rm SU}(2)$ we
  conclude that $\varpi=-2\pi\omega$.
   So, the Hamiltonian groups of $({\mathbb C}P^1,\,\omega)$ and $({\mathbb
   C}P^1,\,\varpi)$ are isomorphic.
  By Theorem \ref{ThmpiHam}, $\sharp\pi_1({\rm
   Ham}({\mathbb C}P^1,\,\varpi))\geq 2$. This result is consistent
   with the fact that  $\pi_1({\rm
   Ham}({\mathbb C}P^1,\,\omega))$ is isomorphic to ${\mathbb Z}/2{\mathbb
   Z}$ \cite[page 52]{lP01}.

\end{document}